\documentclass[12pt]{article}

\usepackage{epsf,amsfonts,hyperref}
\renewcommand{\appendix}[1]{
    \addtocounter{section}{1}
    \setcounter{equation}{0}
    \renewcommand{\thesection}{\Alph{section}}
    \section*{Appendix \thesection\protect\indent #1}
    \addcontentsline{toc}{section}{Appendix \thesection\ \ \ #1}
}
\newcommand\encadremath[1]{\vbox{\hrule\hbox{\vrule\kern8pt
\vbox{\kern8pt \hbox{$\displaystyle #1$}\kern8pt}
\kern8pt\vrule}\hrule}}
\def\enca#1{\vbox{\hrule\hbox{
\vrule\kern8pt\vbox{\kern8pt \hbox{$\displaystyle #1$}
\kern8pt} \kern8pt\vrule}\hrule}}

\newcommand\figureframex[3]{
\begin{figure}[bth]
\hrule\hbox{\vrule\kern8pt
\vbox{\kern8pt \vbox{
\begin{center}
{\mbox{\epsfxsize=#1.truecm\epsfbox{#2}}}
\end{center}
\caption{#3}
}\kern8pt}
\kern8pt\vrule}\hrule
\end{figure}
}
\newcommand\figureframey[3]{
\begin{figure}[bth]
\hrule\hbox{\vrule\kern8pt
\vbox{\kern8pt \vbox{
\begin{center}
{\mbox{\epsfysize=#1.truecm\epsfbox{#2}}}
\end{center}
\caption{#3}
}\kern8pt}
\kern8pt\vrule}\hrule
\end{figure}
}

\renewcommand{\thesection}{\arabic{section}}

\makeatletter
\@addtoreset{equation}{section}
\makeatother
\newtheorem{theorem}{Theorem}[section]

\newtheorem{remark}{Remark}[section]
\newtheorem{proposition}{Proposition}[section]
\newtheorem{lemma}{Lemma}[section]
\newtheorem{corollary}{Corollary}[section]
\newtheorem{definition}{Definition}[section]
\def\br{\begin{remark}\rm\small}
\def\er{\end{remark}}
\def\bt{\begin{theorem}}
\def\et{\end{theorem}}
\def\bd{\begin{definition}}
\def\ed{\end{definition}}
\def\bp{\begin{proposition}}
\def\ep{\end{proposition}}
\def\bl{\begin{lemma}}
\def\el{\end{lemma}}
\def\bc{\begin{corollary}}
\def\ec{\end{corollary}}
\def\beaq{\begin{eqnarray}}
\def\eeaq{\end{eqnarray}}
\newcommand{\proof}[1]{{\noindent \bf proof:}\par
{#1} $\square$
\bigskip}

\newcommand{\beq}{\begin{equation}}
\newcommand{\eeq}{\end{equation}}
\newcommand{\bea}{\begin{eqnarray}}
\newcommand{\eea}{\end{eqnarray}}

\renewcommand{\and}{{\qquad {\rm and} \qquad}}

\newcommand{\virg}{{\qquad , \qquad}}

 \newcommand{\Tr}{{\,\rm Tr}\:}

\newcommand{\Res}{\mathop{\,\rm Res\,}}

\newcommand{\td}[1]{{\tilde{#1}}}

\renewcommand{\l}{\lambda}

\newcommand{\ee}[1]{{{\rm e}^{#1}}}

\newcommand{\R}{{\mathbf R}}

\newcommand{\Pint}{{\int\kern -1.em -\kern-.25em}}

\renewcommand{\l}{\lambda}
\renewcommand{\L}{\Lambda}

\newcommand{\curve}{{\cal E}}

\begin{document}
\sloppy


\pagestyle{empty}
\hfill SPT-07/104
\addtolength{\baselineskip}{0.20\baselineskip}
\begin{center}
\vspace{26pt}
{\large \bf {Recursion between Mumford volumes of moduli spaces}}
\newline
\vspace{26pt}

{\sl B.\ Eynard}\hspace*{0.05cm}\footnote{ E-mail: bertrand.eynard@cea.fr }\\
\vspace{6pt}
Service de Physique Th\'{e}orique de Saclay,\\
F-91191 Gif-sur-Yvette Cedex, France.\\
\end{center}

\vspace{20pt}
\begin{center}
{\bf Abstract}

We propose a new proof, as well as a generalization of Mirzakhani's recursion for volumes of moduli spaces.
We interpret those recursion relations in terms of expectation values in Kontsevich's integral, i.e. we relate them to a Ribbon graph decomposition of Riemann surfaces.
We find a generalization of Mirzakhani's recursions to measures containing all higher Mumford's $\kappa$ classes, and not only $\kappa_1$ as in the Weil-Petersson case.

\end{center}
%





\vspace{26pt}
\pagestyle{plain}
\setcounter{page}{1}


\section{Introduction}

Let
\beq
 {\rm Vol}_{{\rm WP}}({\cal M}_{g,n}(L_1,\dots,L_n))
\eeq
be the volume (measured with Weil-Petersson's measure) of the moduli space of genus $g$ curves with $n$ geodesic boundaries of length $L_1,\dots,L_n$.
Maryam Mirzakhani found a beautiful recursion relation \cite{Mirza1,Mirza2} for those functions, allowing to compute all of them in principle. That relation has then received several proofs \cite{mulsaf, Xu}, and we provide one more proof, more ``matrix model oriented''.

The main interest of our method, is that it easily generalizes to a larger class of measures, containing all Mumford classes $\kappa$, which should also prove the result of Liu and Xu \cite{Xu}. 

In fact, our recursion relations are those of \cite{EOFg}, and they should be generalizable to a much larger set of measures, not only those based on Kontsevich's hyperelliptical spectral curve, and not only rational spectral curves.
For instance they hold for the generalized Kontsevich integral whose spectral curve is not hyperelliptical, i.e. they should hopefully allow to compute also some sort of volumes of moduli spaces of stable maps with spin structures.

\smallskip

In \cite{EOWP} it was observed that after Laplace transform, Mirzakhani's recursion became identical to the solution of loop equations \cite{EOFg} for Kontsevich's matrix integral.
Based on that remark we are in position to reprove Mirzakhani's result, and in fact we prove something more general:

Consider an arbitrary set of Kontsevich KdV times\footnote{Our definition of times $t_k$ slightly differs from the usual one, we have $t_k={1\over N}\Tr \L^{-k}$.} $t_{2d+3}$, $d=0,1,\dots,\infty$, we define their conjugated times $\td{t}_k$, $k=0,1,\dots,\infty$, by:
\beq
f(z)  = \sum_{a=1}^\infty {(2 a+1)!\over a!}\,\,{t_{2a+3}\over 2-t_3} \,\, z^a 
\quad \to \quad
\td{f}(z)= -\ln{(1-f(z))} = \sum_{b=1}^\infty \td{t}_b \,\, z^b  
\eeq
Then we prove the following theorem:

\bt\label{mainth}
Given a set of conjugated Kontsevich times $\td{t}_0,\td{t}_1,\td{t}_2,\dots$, the following ``Mumford volumes'', 
 \bea\label{defWgnmainth}
  W_{g,n}(z_1,\dots,z_n)  
&=& 2^{- d_{g,n}}(t_3-2)^{2-2g-n}\!\!\!\!  \sum_{d_0+d_1+\dots+d_n=d_{g,n}}\, 
 \sum_{k=1}^{d_0} {1\over k!}\,\sum_{b_1+\dots+b_k =d_0, b_i>0}  \cr
 && \qquad \qquad \prod_{i=1}^n {2d_i+1!\over d_i!}\, {dz_i\over z_i^{2d_i+2}}\,\, \prod_{l=1}^k \td{t}_{b_l} <\prod_{l=1}^k \kappa_{b_l} \prod_{i=1}^n  \psi_i^{d_i}>_{g,n}  \cr
\eea
where $d_{g,n}=3g-3+n={\rm dim}\, {\cal M}_{g,n}$,
satisfy the following recursion relations (where $K=\{z_1,\dots,z_n\}$):
\bea\label{recWgninmainth}
&& W_{0,1}=0
\qquad
W_{0,2}(z_1,z_2)={dz_1 dz_2\over (z_1-z_2)^2} \cr
W_{g,n+1}(K,z_{n+1}) 
&=& {1\over 2}\Res_{z\to 0} { dz_{n+1} \over ( z_{n+1}^2-z^2)(y(z)-y(-z))dz}\,\,\Big[ W_{g-1,n+2}(z,-z,K) \cr
&& \qquad\qquad 
 + \sum_{h=0}^g\sum_{J\subset K} W_{h,1+|J|}(z,J)\,W_{g-h,1+n-|J|}(-z,K/J)
\Big] \cr
\eea
where
\beq
y(z) = z-{1\over 2}\sum_{k=0}^\infty t_{2k+3}\,z^{2k+1}
\eeq

\et

From theorem.\ref{mainth}, we obtain as an immediate consequence if  $t_{2d+3} = -{(2i\pi)^{2d}\over 2d+1!} + 2\delta_{d,0}$, i.e. $\td{t}_1=4\pi^2$ and $\td{t}_k=0$ for $k>1$, and after Laplace transform:
\bc
The Weil-Petersson volumes satisfy Mirzakhani's recursions.
\ec

\bigskip

The proof of theorem \ref{mainth} is detailed in the next sections, it can be sketched as follows:
\begin{itemize}
\item We first define some $W_{g,n}(z_1,\dots,z_n)$ which obey the recursion relations of \cite{EOFg}, i.e. eq.\ref{recWgninmainth}. In other words, we define them as the solution of the recursion, without knowing what they compute.

\item We prove that those $W_{g,n}(z_1,\dots,z_n)$ correspond to some expectation values in the Kontsevich integral $Z(\L) = \int dM \,\, \ee{-N\Tr({M^3\over 3} - M \L^2)}$, where $\L={\rm diag}(\l_1,\dots,\l_n)$, and $t_k={1\over N}\Tr \L^{-k}$, of the form:
\beq
{W_{g,n} (\l_{i_1},\dots,\l_{i_n}) \over dx(\l_{i_1})\dots dx(\l_{i_n})} = (-1)^n \left< M_{i_1,i_1}\dots M_{i_n,i_n}\right>_c^{(g)}
\eeq

\item Then we expand $< M_{i_1,i_1}\dots M_{i_n,i_n}>$ into Feynman ribbon graphs, which are in bijection with a cell decomposition of ${\cal M}_{g,n}^{\rm comb}$ (like in Kontsevich's first works), and the value of each of those Feynman graphs is precisely the Laplace transform of the volume of the corresponding cell.

\item the sum over all cells yields the expected result: the inverse Laplace transforms of $W_{g,n}$ are the volumes $V_{g,n}$, and, by definition, they satisfy the recursion relations.

\item In fact the volumes are first written in terms of the first Chern classes $\psi_i$ in formula eq.\ref{eqVgnintersec}, and after some combinatorics, we find more convenient to rewrite them in terms of Mumford $\kappa$ classes.

\end{itemize}

\bigskip

Then, we specialize our theorem to some choices of times $t_k$'s, in particular the following:

$\bullet$ The first example is $t_{2d+3} = -{(2i\pi)^{2d}\over 2d+1!} + 2\delta_{d,0}$, in which case $V_{g,n}$ the Laplace transform of $W_{g,n}$ are the Weil-Petersson volumes, and thus we recover Mirzakhani's recursions.

$\bullet$ Our second example is $t_k=\l^{-k}$, i.e. $\L=\l\, {\rm Id}$, for which the Kontsevich integral reduces to a standard one-matrix model, and for which the $W_{g,n}$ are known to count triangulated maps, i.e. discrete surfaces with the discrete Regge metrics (metrics whose curvature is localized on vertices and edges).
We are thus able to associate some class to that discrete measure on ${\cal M}_{g,n}$.
And we have a formula which interpolates between the enumeration of maps and the enumeration of Riemann surfaces, in agreement with the spirit of 2d-quantum gravity in the 80's \cite{witten, BIPZ, DGZ}.

\section{Proof of the theorem}

\subsection{Kontsevich's integral}

In his very famous work \cite{Konts} Maxim Kontsevich introduced the following matrix integral as a generating function for intersection numbers
\bea\label{defZ}
Z(\L) &=&\int dM \,\, \ee{-N\Tr({M^3\over 3} - M(\L^2+t_1))} \cr
   &=&\ee{{2N\over 3}\Tr \L^3+N t_1\Tr\L}\,\,\int dM \,\, \ee{-N\Tr({M^3\over 3} + M^2\L  - t_1 M  )} \cr
\eea
where the integral is a formal integral over hermitian matrices $M$ of size $N$, and  $\L$ is a fixed diagonal matrix
\beq
\L={\rm diag}(\l_1,\dots,\l_n)
\virg
t_k={1\over N}\Tr \L^{-k}
\eeq
Throuthough all this article we shall assume $t_1=0$, since anyways none of the quantities we are interested in here depend on $t_1$ (see symplectic invariance in \cite{EOFg}, or see \cite{DFIZ}).

\medskip

In \cite{EOFg}, a method to compute the topological expansion of such matrix integrals was developped.
We first define the Kontsevich's spectral curve:

\bd
The spectral curve of $Z(\L)$ is the rational plane curve of equation:
\beq
\curve(x,y) = y^2-x -  {y\over N}\,\Tr{1\over x-t_1-\L^2} - {1\over N}\,\left<\Tr{1\over x-t_1-\L^2}\,M\right>^{(0)} = 0
\eeq
i.e. it has the following rational uniformization
\beq
\curve(x,y) =
\left\{\begin{array}{l}
x(z) = z^2+t_1
\cr
y(z) = z + {1\over 2N}\,\Tr{1\over \L(z-\L)} = z - {1\over 2}\sum_{k=0}^\infty t_{k+2} z^k
\end{array}\right.
\eeq
\ed

Then we define (i.e. the algebraic invariants of \cite{EOFg}):
\bd\label{defrec} We define the correlators:
\beq
W_{0,1}=0
\qquad
W_{0,2}(z_1,z_2)={dz_1 dz_2\over (z_1-z_2)^2}
\eeq
and we define by recursion on $2g-2+n$, the symmetric\footnote{The non-obvious fact that this is symmetric in its $n+1$ variables is proved by recursion in \cite{EOFg}.} form $W_{g,n+1}(z_0,z_1,\dots,z_n)$ by (we write $K=\{z_1,\dots,z_n\}$): 
\bea
W_{g,n+1}(K,z_{n+1}) 
&=& \Res_{z\to 0} {z\, dz_{n+1} \over ( z_{n+1}^2-z^2)(y(z)-y(-z))dx(z)}\,\,\Big[ W_{g-1,n+2}(z,-z,K) \cr
&& \qquad\qquad 
 + \sum_{h=0}^g\sum_{J\subset K} W_{h,1+|J|}(z,J)\,W_{g-h,1+n-|J|}(-z,K/J)
\Big] \cr
\eea
Then, if $d\Phi=y dx$, we define for $g>1$:
\beq
F_g = {1\over 2g-2}\,\, \Res_{z\to 0} \Phi(z) W_{g,1}(z) 
\eeq
(there is a separate definition of $F_g$ for $g=0,1$, but we shall not use it here).
\ed

We recall the result of \cite{EOFg} (which uses also \cite{SymFg}):
\bt
\beq
\ln{Z} = \sum_{g=0}^\infty N^{2-2g} F_g
\eeq
\et

Now, we prove the more elaborate result:
\bt
if $i_1,\dots,i_n$ are $n$ distinct integers in $[1,N]$, then:
\beq
\encadremath{
{W_n^{(g)} (\l_{i_1},\dots,\l_{i_n}) \over dx(\l_{i_1})\dots dx(\l_{i_n})} =  \left< M_{i_1,i_1}\dots M_{i_n,i_n}\right>_c^{(g)}
}\eeq
\et
where $<.>$ means the formal expectation value with respect to the measure used to define $Z$, 
the subscript $c$ means connected part or cumulant, and the subscript $(g)$ means the $g^{\rm th}$ term in the $1/N^2$ topological expansion.

In other words, the $W_{g,n}$ compute some expectation values in the Kontsevich integral, which are not the same as those computed by \cite{DFIZ}.

\smallskip
\proof{


From eq. \ref{defZ}, it is easy to see that:
\beq
N^{-n}\,{\partial^n \ln{Z} \over \partial \l_{i_1} \dots \partial \l_{i_n}} =   2^n\, \l_{i_1}\dots\l_{i_n}\,\, \left< M_{i_1,i_1}\dots M_{i_n,i_n} \right>_c
\eeq
i.e., to order $N^{2-2g-n}$:
\beq
{\partial^n F_g \over \partial \l_{i_1} \dots \partial \l_{i_n}} =   2^n\, \l_{i_1}\dots\l_{i_n}\,\, \left< M_{i_1,i_1}\dots M_{i_n,i_n} \right>_c^{(g)}
\eeq
Now, let us compute ${\partial F_g \over \partial \l_i}$ with the method of \cite{EOFg}.

Consider an infinitesimal variation of the matrix $\L$: $\l_i\to \l_i+\delta \l_i$ (we assume $\delta t_1=0$).
It translates into the following variations of the function $y(z)$:
\beq
\delta y(z) = {1\over 2Nz} \,\Tr{\delta \L\over (z-\L)^2}
\eeq
and thus the form:
\bea
- \delta y(z) dx(z) 
&=& d\left({1\over N} \,\Tr{\delta \L\over z-\L}\right) \cr
&=& \Res_{\zeta\to z} {1\over (z-\zeta)^2}\,{1\over N} \,\Tr{\delta \L\over \zeta-\L} \cr
&=& - \sum_i \Res_{\zeta\to \l_i} {1\over (z-\zeta)^2}\,{1\over N} \,\Tr{\delta \L\over \zeta-\L} \cr
\eea
Then, using theorem 5.1 of \cite{EOFg}, we have:
\bea
\delta F_g 
&=&  \sum_i \Res_{\zeta\to \l_i} W_1^{(g)}(\zeta)\,{1\over N} \,\Tr{\delta \L\over \zeta-\L}  \cr
&=& \sum_i   {W_1^{(g)}(\l_i)\over d\l_i}\,{\delta\l_i\over N}   
\eea
i.e.
\beq
W_1^{(g)}(\l_i) =  \left< M_{ii} \right>^{(g)}\,dx(\l_i) 
\eeq
And repeating the use of theorem 5.1 in \cite{EOFg} recursively we get the result.
}

Example:
\beq
\left< M_{ii} \right>^{(1)} = {1\over 16(2-t_3)}\,\left({1\over \l_i^5} + {t_5\over (2-t_3)\l_i^3}\right)
\quad \longrightarrow \quad
\left< \Tr M \right>^{(1)} = {t_5\over 8(2-t_3)^2}
\eeq

\subsection{Expectation values and ribbon graphs}

Let $i_1,\dots,i_n$ be $n$ distinct given integers $\in [1,\dots,N]$. We want to compute:
\beq
\left< M_{i_1,i_1} \dots M_{i_n,i_n} \right>^{(g)}
\eeq
Let us also choose $n$ positive real perimeters $P_1,\dots,P_n$

Let $\Gamma(g,n,m)$ be the set of tri-valent oriented ribbon graphs of genus $g$, with $n$ marked faces, and $m$ unmarked faces. Each marked face $F=1,\dots,n$ carries the given index $i_F$, and each unmarked face $f$ carries an index $i_f\in [1,\dots,N]$.

Let us consider another set of graphs:
Let $\Gamma^*(g,n,m)$ be the set of oriented ribbon graphs of genus $g$, with trivalent and 1-valent vertices, made of $m$ unmarked faces bordered with only tri-valent vertices, each of them carrying an index $i_f$,  
and $n$ marked faces carrying the fixed index $i_F\in\{ i_1,\dots,i_n\}$, such that each marked face has one 1-valent vertex on its boundary.
The unique trivalent vertex linked to the 1-valent vertex on each marked face, corresponds to a marked point on the boundary of that face.

$$
{\mbox{\epsfxsize=4.truecm\epsfbox{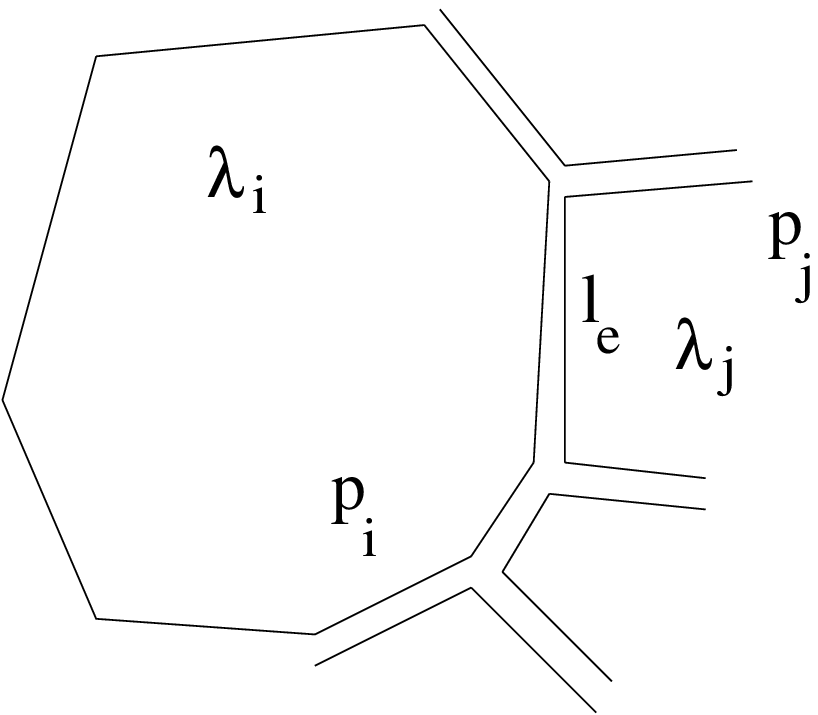}}}
\qquad , \qquad
{\mbox{\epsfxsize=4.truecm\epsfbox{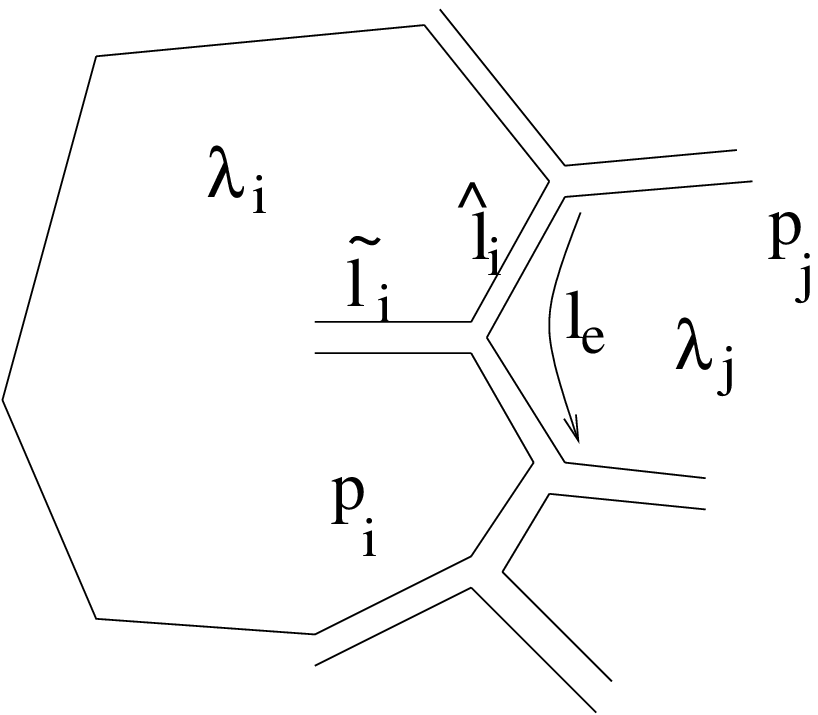}}}
$$


For any graph $G$ in either $\Gamma(g,n,m)$ or $\Gamma^*(g,n,m)$, each edge $e$ is bordered by two faces (possibly not different), and we denote the pair of their indices as $(e_{\rm left},e_{\rm right})$.

\bigskip

Assume that $i_1,\dots,i_n$ are distinct integers.
The usual fat graph expansion of matrix integrals gives (cf \cite{BIPZ, DGZ, Konts}):
\beq
\left< M_{i_1,i_1} \dots M_{i_n,i_n} \right>^{(g)} 
= N^{-m}\,\sum_m \sum_{G\in \Gamma^*_{g,n,m}}  \sum_{\{i_f\} }\, {(-1)^{\# {\rm vertices}}\over \#{\rm Aut}(G)}\,\,\prod_{e\in {\rm edges}(G)}\, {1\over \l_{e\,{\rm left}}+\l_{e\,{\rm right}}} 
\eeq
It is obtained by first expanding $\ee{-{N\over 3}\Tr M^3}=\sum_{v=0}^\infty {N^v\over 3^v\, v!} \,\, (-1)^v \,\, (\Tr M^3)^v$, and then computing each polynomial moment of the Gaussian measure $\ee{-N\Tr \L M^2}$ with the help of Wick's theorem.
Each $\Tr M^3$ corresponds to a trivalent vertex, each $M_{ii}$ corresponds to a 1-valent vertex, and edges correspond to the ``propagator'' $<M_{ij}M_{kl}>_{{\rm Gauss}} = {\delta_{il}\delta_{jk}\over N(\l_i+\l_j)}$. 
The result is best represented as a fatgraph, whose edges are double lines, carrying two indices.
The indices are conserved along simple lines.
The symmetry factor comes from the combination of $1/(3^v v!)$ and the fact that some graphs are obtained several times.
Notice that $(-1)^v=(-1)^n$, because the total number of 1 and 3-valent vertices must be even.

\medskip

Notice that the edge connected to the 1-valent vertex $M_{i_F,i_F}$ gives a factor $1/2\l_{i_F}$, and the two edges on the boundary of face $F$, on each side of the 1-valent vertex give a factor $1/(\l_{i_F}+\l_{j})^2$ (where $j$ is the index of the neighboring face), which can be written:
\beq
{1\over (\l_{i_F}+\l_j)^2} = \int_0^\infty dl_e \int_0^{l_e} d\hat{l}_i\,\, \ee{-l_e(\l_{i_F}+\l_j)}
\eeq
and all other edges have a weight of the form:
\beq
{1\over \l_{e\, {\rm left}}+\l_{e\, {\rm right}}} = \int_0^\infty dl_e \,\, \ee{-l_e(\l_{e\, {\rm left}}+\l_{e\, {\rm right}})}
\eeq
We are thus led to associate to each edge $e$ a length $l_e\in {\bf R}^+$.

Therefore
\bea
 \left< M_{i_1,i_1} \dots M_{i_n,i_n} \right>^{(g)} 
&=& {N^{-m}\,\over 2^n \l_{i_1}\dots \l_{i_n}}\,\sum_m \sum_{G\in \Gamma^*_{g,n,m}}  \sum_{\{i_f\} }\,{(-1)^n\over \#{\rm Aut}(G)} \,\, \cr
&& \qquad \prod_{e\in {\rm edges}(G)}\,\, \int_0^\infty dl_e\,\, \ee{- \sum_e l_e(\l_{e\,{\rm left}}+\l_{e\,{\rm right}})} \,\,\,\prod_{F=1}^n\int_0^{l_{e_F}} d\hat{l}_F\cr
\eea
Now, we introduce the perimeters of each face $P_F$ for marked faces, and $p_f$ for unmarked ones.

Notice that each graph of $\Gamma^*_{g,n,m}$ projects on a graph of $\Gamma_{g,n,m}$ by removing the 1-valent vertex and its adjacent trivalent vertex, and keeping a marked point on the boundary of the face $F$. 
The sum of $\int \prod_F d\hat{l}_F$ over graphs of $\Gamma^*_{g,n,m}$ which project to the same graph, corresponds to a sum of all possibilities of marking a point on the boundary of face $F$, i.e. a factor $P_F$, and thus removing the marked point.
Therefore:
\bea
&& \left< M_{i_1,i_1} \dots M_{i_n,i_n} \right>^{(g)} \cr
&=& {N^{-m}\,\over 2^n \l_{i_1}\dots \l_{i_n}}\,\sum_m \sum_{G\in \Gamma_{g,n,m}}  \sum_{\{i_f\} }\,{(-1)^n\over \#{\rm Aut}(G)} \,\,  \prod_f \int_0^\infty dp_f \,\, \ee{-\sum_{f} \l_{i_f} p_f}  \cr
&& \,\,  \prod_F \int_0^\infty P_F \, dP_F\,   \ee{-\l_{i_F}P_F}
\,\, \prod_e \int_0^\infty  \,dl_e\,\,   \prod_{f} \delta(p_f-\sum_{e\in \partial f} l_e) \,\prod_{F=1}^n \delta(P_F-\sum_{e\in \partial F} l_e) \cr
&=& {1\over 2^n \l_{i_1}\dots \l_{i_n}}\,\sum_m \sum_{G\in \Gamma_{g,n,m}} {(-1)^n\over \#{\rm Aut}(G)} \,\,  \cr
&& \prod_f \int_0^\infty dp_f \,\,  {1\over N}\Tr(\ee{- p_f \L})  \prod_F \int_0^\infty P_F\, dP_F\,  \ee{-\l_{i_F}P_F} \,\,\,\, {\rm Vol}(\pi_G^{-1}(P_F,p_f)) \cr
\eea
where ${\rm Vol}(\pi_G^{-1}(P_F,p_f))$ is the volume of the pullback of the ribbon graph $G$ in ${\cal M}_{g,n+m}^{\rm comb}$:
\beq
{\rm Vol}(\pi_G^{-1}(P_F,p_f)) = \int \prod_e   \,dl_e\,\,   \prod_{f} \delta(p_f-\sum_{e\in \partial f} l_e) \,\prod_{F=1}^n \delta(P_F-\sum_{e\in \partial F} l_e)
\eeq
The number of integrations (i.e. after performing the $\delta$) is
$2d_{g,n+m}=\#{\rm edges}-\#{\rm faces}=2(3g-3+n+m)$, which is the dimension of ${\cal M}_{g,n+m}$, therefore $\prod_e dl_e$ is a top-dimension volume form on ${\cal M}_{g,n+m}^{\rm comb}=\overline{\cal M}_{g,n+m}\times \R_+^{n+m}$, i.e.:
\beq
\prod_e   \,dl_e = {\rho_{g,n+m}\over d_{g,n+m}!}\,\prod_F dP_F \prod_f dp_f
 \,\,\,\wedge  \Omega^{d_{g,n+m}}
 \eeq
where $\Omega$ is the 2-form on the strata $\pi^{-1}_G(P_F,p_f)$ of ${\cal M}_{g,n+m}^{\rm comb}$ such that:
\beq
\Omega = \sum_f p_f^2 \omega_f + \sum_F P_F^2 \omega_F 
\eeq
and where $\omega_f = \sum_{e<e'} d(l_e/p_f)\wedge d(l_{e'}/p_f)$ is the first Chern class of pullback of the cotangent bundle at the center of the face  $\psi_f=c_1({\cal L}_f)$.

Kontsevich \cite{Konts} proved that the constant $\rho_{g,n+m}$ is given by:
\beq
\rho_{g,n+m}= 2^{g-1 - 2d_{g,n+m}}
\eeq

Thus we have:
\bea
&&  {\rm Vol}(\pi^{-1}_G(P_F,p_f))  \cr
&=& {\rho_{g,n+m}\over d_{g,n+m}!}\,\int_{\pi^{-1}_G(P_F,p_f)} \Omega^{d_{g,n+m}}  \cr
&=& \rho_{g,n+m}\, \sum_{\sum_f d_f + \sum_F d_F=d_{g,n+m}}\,\int_{\pi^{-1}_G(P_F,p_f)} \, \prod_f {p_f^{2 d_f}\,\psi_f^{d_f}\over d_f !}\,\prod_F {P_F^{2 d_F}\,\psi_F^{d_F}\over d_F !}  \cr
&=& \rho_{g,n+m}\, \sum_{\sum_f d_f + \sum_F d_F=d_{g,n+m}}\, \, \prod_f {p_f^{2 d_f}\over d_f !}\,\prod_F {P_F^{2 d_F}\over d_F !}\,\,<\prod_f \psi_f^{d_f} \prod_F  \psi_F^{d_F}>_G  \cr
\eea
therefore
\bea
&& \prod_f \int_0^\infty dp_f \,\,  {1\over N}\Tr(\ee{- p_f \L}) {\rm Vol}(\pi^{-1}_G(P_F,p_f))  \cr
&=& \rho_{g,n+m}\, \,\prod_f \int_0^\infty dp_f \,\,  {1\over N}\Tr(\ee{- p_f \L}) \sum_{\sum_f d_f + \sum_F d_F=d_{g,n+m}}\cr
&& \quad \, \prod_f {p_f^{2 d_f}\over d_f !}\,\prod_F {P_F^{2 d_F}\over d_F !}\,\,<\prod_f \psi_f^{d_f} \prod_F  \psi_F^{d_F}> _G   \cr
&=& \rho_{g,n+m}\,\,   \sum_{\sum_f d_f +\sum_F d_F=d_{g,n+m}} \cr
&& \quad \, \, \prod_f {2d_f!\over d_f !}\,\,{1\over N}\Tr(\L^{-(2d_f+1)}) \,\prod_F {P_F^{2 d_F}\over d_F !}\,\,<\prod_f \psi_f^{d_f} \prod_F  \psi_F^{d_F}>_G    \cr
&=& \rho_{g,n+m}\,\,   \sum_{\sum_f d_f + \sum_F d_F=d_{g,n+m}}\, \, \prod_f {2d_f!\over d_f !}\,\,t_{2d_f+1} \,\prod_F {P_F^{2 d_F}\over d_F !}\,\,<\prod_f \psi_f^{d_f} \prod_F  \psi_F^{d_F}>_G    \cr
\eea

and then, when we sum over all graphs (since we sum over graphs with $m$ unmarked faces, we have to divide wrt to the symmetry factor $m!$, like in \cite{Konts}) :
\bea
&& \left< M_{i_1,i_1} \dots M_{i_n,i_n} \right>^{(g)} \cr
&=& {(-1)^n \rho_{g,n}\over 2^n \l_{i_1}\dots \l_{i_n}}\, \prod_F \int_0^\infty P_F dP_F  \ee{-\l_{i_F}P_F} \,\, \cr
&& \qquad \sum_m  {1\over m!}\,\sum_{\sum d_f + d_F=d_{g,n}+m}\,   \,\,   \prod_f {2d_f!\over d_f !}\,\,{t_{2d_f+1}\over 4} \,\prod_F {P_F^{2 d_F}\over d_F !}\,\,<\prod_f \psi_f^{d_f} \prod_F  \psi_F^{d_F}> \cr
\eea

Therefore, if we write:
\bea
{W_{g,n}(\l_{i_1},\dots,\l_{i_n}) \over d\l_{i_1}\dots d\l_{i_n}} = \int_0^\infty dP_1 \dots dP_n \,\,\, \prod_F P_F\,\ee{-\l_{i_F} P_F} \,\, V_{g,n}(P_1,\dots,P_n)
\eea
we find that the inverse Laplace transform of $W_{g,n}$ is:
\beq\label{eqVgnintersec}
\encadremath{
\begin{array}{l}
 V_{g,n}(P_1,\dots,P_n) 
= \rho_{g,n}\, \, \sum_m  {(-1)^n\over m!}\,\sum_{\sum_1^m d_f + \sum_1^n d_F=d_{g,n}+m} \cr
\,   \,\,   \prod_f {2d_f!\over d_f !}\,\,{t_{2d_f+1}\over 4} \,\prod_F {P_F^{2 d_F}\over d_F !}\,\,<\prod_f \psi_f^{d_f} \prod_F  \psi_F^{d_F}> 
\end{array}
}\eeq
where the intersection theory is computed on $\overline{\cal M}_{g,n+m}$.

Since we are interested only in the perimeters of the $n$ marked faces, we may try to perform the integration over the $m$ unmarked faces, i.e. we introduce the forgetful projection $\pi_{n+m\to n}:\overline{\cal M}_{g,n+m}\to \overline{\cal M}_{g,n}$ which ``forgets'' the $m$ remaining points. 
It is known \cite{arbcor, witten} that the push forward of the classes $\psi_f^{d_f}$, can then be rewritten in terms of Mumford's \cite{mumford} tautological classes $\kappa_b$ on $\overline{\cal M}_{g,n}$, by the relation:

\beq
(\pi_{n+m\to n})_*(\psi_1^{a_1+1}\dots\psi_m^{a_m+1} \prod_F \psi_F^{d_F}) = \sum_{\sigma\in \Sigma_m}\,\, \prod_{c={\rm cycles\, of\,}\sigma} \,\,\, \kappa_{\sum_{i\in c} a_i}\,\prod_F \psi_F^{d_F}
\eeq

Therefore, if we rewrite $d_f=a_f+1$ we have:
\bea
&&  {1\over \rho_{g,n}}\, V_{g,n}(P_1,\dots,P_n)  \cr
&=&  \sum_m  {(-1)^n\over m!}\,\sum_{\sum_1^m a_f + \sum_1^n d_F=d_{g,n}} \cr
&& \qquad \,   \,\,   \prod_f {2 a_f+1!\over a_f !}\,\,{t_{2a_f+3}\over 2} \,\prod_F {P_F^{2 d_F}\over d_F !}\,\,<\prod_f \psi_f^{a_f+1} \prod_F  \psi_F^{d_F}>  \cr
&=&  \sum_m  {(-1)^n\over m!}\,\sum_{\sum_1^m a_f + \sum_1^n d_F=d_{g,n}}\,  \sum_{\sigma\in\Sigma_m} \cr
&& \qquad \,\,   \prod_f {2 a_f+1!\over a_f !}\,\,{t_{2a_f+3}\over 2} \,\prod_F {P_F^{2 d_F}\over d_F !}\,\,<\prod_c \kappa_{\sum_c a_i} \prod_F  \psi_F^{d_F}>  \cr
&=&  (-1)^n\,\sum_{d_0+d_1+\dots+d_F=d_{g,n}}\, \prod_F {P_F^{2 d_F}\over d_F !}\,\,  \sum_m  {1\over m!}\,\sum_{a_1+\dots+a_m =d_0, a_f\geq 0}\,  \sum_{\sigma\in\Sigma_m} \cr
&& \qquad \,\,   \prod_f {2 a_f+1!\over a_f !}\,\,{t_{2a_f+3}\over 2} \,\,<\prod_c \kappa_{\sum_c a_i} \prod_F  \psi_F^{d_F}>  
\eea
Now, instead of summing over permutations, let us sum over classes of permutations, i.e. partitions $l_1\geq l_2\geq \dots \geq l_k >0$, and we denote $|l|=\sum_i l_i=m$ the weight of the class, and $|[l]|$ the size of the class:
\beq
|[l]| = {|l|! \over \prod_i l_i\,\, \prod_j (\# \{i/\,\, l_i=j\})!}
\eeq
The sum over the $a$'s for each class gives:
\bea
&&  {(-1)^n\over \rho_{g,n}}\, V_{g,n}(P_1,\dots,P_n)  \cr
&=&  \sum_{d_0+d_1+\dots+d_F=d_{g,n}}\, \prod_F {P_F^{2 d_F}\over d_F !}\,\, \sum_k \, \sum_{l_1\geq l_2\geq \dots \geq l_k >0}\, \, {|[l]|\over |l|!} \, \,\sum_{a_{i,j}, i=1,\dots,k, j=1,\dots,l_i} \cr
&& \qquad \,\, \delta(\sum_{i,j} a_{i,j}-d_0) \,\,   \prod_{i,j} {2 a_{i,j}+1!\over a_{i,j} !}\,\,{t_{2a_{i,j}+3}\over 2} \,\,<\prod_{i=1}^k \kappa_{\sum_{j=1}^{l_i} a_{i,j}} \prod_F  \psi_F^{d_F}>  \cr
\eea
Since the summand is symmetric in the $l_i$'s, the ordered sum over $l_1\geq \dots l_k$, can be replaced by an unordered sum (multiplying by $1/k!$, and by $\prod_i (\# \{i/\,\, l_i=j\})!$ in case some $l_i$ coincide):
\bea
&&  {(-1)^n\over \rho_{g,n}}\, V_{g,n}(P_1,\dots,P_n)  \cr
&=&  \sum_{d_0+d_1+\dots+d_F=d_{g,n}}\, \prod_F {P_F^{2 d_F}\over d_F !}\,\, \sum_k {1\over k!}\, \sum_{l_1,l_2,\dots,l_k >0}\, \, \prod_{i=1}^k {1\over l_i} \,\, \,\sum_{a_{i,j}, i=1,\dots,k, j=1,\dots,l_i}\cr
&& \qquad \,\, \delta(\sum_{i,j} a_{i,j}-d_0) \,\,   \prod_{i,j} {2 a_{i,j}+1!\over a_{i,j} !}\,\,{t_{2a_{i,j}+3}\over 2} \,\,<\prod_{i=1}^k \kappa_{\sum_{j=1}^{l_i} a_{i,j}} \prod_F  \psi_F^{d_F}>  \cr
&=&  \sum_{d_0+d_1+\dots+d_F=d_{g,n}}\, \prod_F {P_F^{2 d_F}\over d_F !}\,\, \sum_k {1\over k!}\, \sum_{b_1+b_2+\dots+b_k=d_0}\, \, \prod_{i=1}^k \td{t}_{b_i}  \,\,<\prod_{i=1}^k \kappa_{b_i} \prod_F  \psi_F^{d_F}>  \cr
\eea
where
\beq
\td{t}_b  = \sum_{l>0}\,  {1\over l}  \,\sum_{a_1+\dots+a_{l}=b}\,\,  \prod_{j} {2 a_{j}+1!\over a_{j} !}\,\,{t_{2a_{j}+3}\over 2}   
\eeq
$\td{t}_b$ can be computed as follows: introduce the generating function
\beq
g(z) = \sum_{a=0}^\infty {2 a+1!\over a !}\,\,{t_{2a+3}\over 2} \,z^a
\eeq
then $\td{t}_b$ is
\beq
\td{t}_b = \sum_{l>0}\,  {1\over l} (g^l)_b  = (-\ln{(1-g)})_b 
\eeq
where the subscript $b$ means the coefficient of $z^b$ in the small $z$ Taylor expansion of the corresponding function, i.e.
\beq
-\ln{(1-g(z))} = \sum_{b=0}^\infty \td{t}_b\, z^b = \td{g}(z)
\virg
1-g(z) = \ee{-\td{g}(z)}
\eeq
In fact, it is better to treat the $a=0$ and $b=0$ terms separately.
Define:
\beq
f(z) = 1- {1-g(z)\over 1-{t_3\over 2}} =  \sum_{a=1}^\infty {2 a+1!\over a !}\,\,{t_{2a+3}\over 2-t_3} \,z^a
\eeq
and
\beq
\td{f}(z) = -\ln{(1-f(z))} = \td{g}(z)-\td{t}_0 = \sum_{b=1}^\infty \td{t}_b\, z^b
\eeq

We have:
\beq
\td{t}_0 = -\ln{(1-{t_3\over2})}
\eeq
and $\td{t}_b$ is now a finite sum:
\beq\label{formulafortdtb}
\td{t}_b  = \sum_{l=1}^b\,  {(-1)^l\over l}  \,\sum_{a_1+\dots+a_{l}=b, a_i>0}\,\,  \prod_{j} {2 a_{j}+1!\over a_{j} !}\,\,{t_{2a_{j}+3}\over t_3-2}   
\eeq

\medskip
 
Using that $\kappa_0 = 2g-2+n$, we may also perform the sum over all vanishing $b$'s.
Let us change $k\to k+l$ where $l$ is the number of vanishing $b$'s, i.e.
\bea  
 {(-1)^n\over \rho_{g,n}}\, V_{g,n}(P_1,\dots,P_n)  
&=& \!\!\!\!  \sum_{d_0+d_1+\dots+d_F=d_{g,n}}\, \prod_F {P_F^{2 d_F}\over d_F !}\,\, \sum_k \sum_l {1\over k! l!} \,(\td{t}_0 \kappa_0)^l \cr
&& \qquad \, \sum_{b_1+b_2+\dots+b_k=d_0, b_i>0}\, \, \prod_{i=1}^k \td{t}_{b_i}  \,\,<\prod_{i=1}^k \kappa_{b_i} \prod_F  \psi_F^{d_F}>  \cr
&=& \!\!\!\!  \ee{\td{t}_0 \kappa_0}\, \sum_{d_0+d_1+\dots+d_F=d_{g,n}}\, \prod_F {P_F^{2 d_F}\over d_F !}\,\, \sum_k {1\over k!} \cr
&& \qquad   \, \sum_{b_1+b_2+\dots+b_k=d_0, b_i>0}\, \, \prod_{i=1}^k \td{t}_{b_i}  \,\,<\prod_{i=1}^k \kappa_{b_i} \prod_F  \psi_F^{d_F}>  \cr
&=& \!\!\!\! \left(2\over 2-t_3\right)^{2g-2+n}\,\, \sum_{d_0+d_1+\dots+d_F=d_{g,n}}\, \prod_F {P_F^{2 d_F}\over d_F !}\,\, \sum_k {1\over k!}  \cr
&& \qquad \, \sum_{b_1+b_2+\dots+b_k=d_0, b_i>0}\, \, \prod_{i=1}^k \td{t}_{b_i}  \,\,<\prod_{i=1}^k \kappa_{b_i} \prod_F  \psi_F^{d_F}>  \cr
\eea

Notice that:
\beq
\rho_{g,n} 2^{2g-2+n} = 2^{-d_{g,n}}
\eeq
thus
\bea  
&& 2^{d_{g,n}}\,(t_3-2)^{2g-2+n}\,\, V_{g,n}(P_1,\dots,P_n)  \cr
&=& \!\!\!\! \!\!\!\!   \sum_{d_0+d_1+\dots+d_F=d_{g,n}}\, \prod_F {P_F^{2 d_F}\over d_F !}\,\, \sum_k {1\over k!}  \, \sum_{b_1+b_2+\dots+b_k=d_0, b_i>0}\, \, \prod_{i=1}^k \td{t}_{b_i}  \,\,<\prod_{i=1}^k \kappa_{b_i} \prod_F  \psi_F^{d_F}>  \cr
\eea

Finaly we obtain theorem \ref{mainth} $\square$.


\section{Examples}

\subsection{Some examples}

First, we give a few examples with general times $t_k$'s.

Using formula \ref{formulafortdtb}, we have:
\beq
\td{t}_1  = -  6\,\,{t_{5}\over t_3-2}   
\quad , \quad
\td{t}_2  = - 60\,\,{t_7\over t_3-2} + 18\,{t_5^2\over (t_3-2)^2}
\eeq
\beq
\td{t}_3 =
- {7!\over 3!}\,{t_9\over t_3-2} + {3! 5!\over 2!}\,{t_5 t_7\over (t_3-2)^2} - {3!^3\over 3}   {t_5^3\over (t_3-2)^3}
\quad ,\, \dots
\eeq

Then we use theorem \ref{mainth} for some examples. In the examples that follow, the first expression is the definition eq.\ref{defWgnmainth}, while the second expression results from the recursion eq.\ref{recWgninmainth}.
\beq
W_{0,3}(z_1,z_2,z_3) 
= {1\over t_3-2}\,{dz_1\,dz_2\,dz_3\over z_1^2\,z_2^2\,z_3^2}<1>_0
 = {1\over t_3-2}\,{dz_1\,dz_2\,dz_3\over z_1^2\,z_2^2\,z_3^2}
\eeq
i.e.
\beq
V_{0,3}(L_1,L_2,L_3) = {1\over t_3-2}
\virg
<1>_0=1
\eeq

\beq
W_{1,1}(z) = {dz\over 2(t_3-2)}\,\left({6\over z^4}<\psi>_1 + {\td{t}_1\over z^2} <\kappa_1>_1 \right) 
= {dz\over 8(t_3-2)}\,\left( {1\over z^4}-{t_5\over (t_3-2)z^2}\right)
\eeq
i.e.
\beq
<\psi>_1={1\over 24}
\virg
<\kappa_1>_1={1\over 24}
\eeq

\bea
W_{1,2}(z_1,z_2) 
&=& {dz_1 dz_2\over 4(t_3-2)^2\, z_1^6 \,z_2^6}\,\Big[ {5!\over 2!} (z_1^4 <\psi_2^2>+z_2^4<\psi_1^2>)+ 3!^2 z_1^2 z_2^2 <\psi_1 \psi_2> \cr
&& \quad + \td{t}_1 z_1^2 z_2^4 <\kappa_1 \psi_1> + \td{t}_1 z_1^4 z_2^2 <\kappa_1 \psi_2>
+ {1\over 2} \td{t}_1^2 z_1^4 z_2^4 <\kappa_1^2 >
+ \td{t}_2 z_1^4 z_2^4 <\kappa_2>
\Big]
\cr
&=& {dz_1 dz_2 \over 8(t_3-2)^4 z_1^6 z_2^6}\,
\Big[ (t_3-2)^2 (5 z_1^4 + 5 z_2^4 + 3 z_1^2 z_2^2) + 6 t_5^2 z_1^4 z_2^4 \cr
&& - (t_3-2)(6 t_5 z_1^4 z_2^2 + 6 t_5 z_1^2 z_2^4 + 5 t_7 z_1^4 z_2^4) \Big]
\eea
i.e.
\beq
<\kappa_1 \psi_1>_1={1\over 2} \virg <\kappa_1^2>_1={1\over 8} \virg <\kappa_2>_1={1\over 24}
\eeq

The recursion equation \ref{recWgninmainth} also gives:
\bea
W_{2,1}(z) &=& - { dz \over 128 (2-t_3)^7 z^{10}} \Big[ 252\, t_5^4 z^8 + 12\, t_5^2 z^6 (2-t_3) (50\, t_7 z^2 + 21\, t_5) \cr
&& \quad + z^4 (2-t_3)^2 ( 252\, t_5^2 + 348\, t_5 t_7 z^2 + 145\, t_7^2 z^4 + 308\, t_5 t_9 z^4) \cr
&& \qquad + z^2 (2-t_3) (203\, t_5 +145\, z^2 t_7 + 105\, z^4 t_9 +105\, z^6 t_{11}) \cr
&& \qquad \quad + 105\, (2 -t_3)^4 \Big] .\cr
\eea

\beq
W_{4,0}(z_1,z_2,z_3,z_4) = 12\,{dz_1 dz_2 dz_3 dz_4\over (t_3-2)^3\, z_1^2 z_2^2 z_3^2 z_4^2}\,
\Big(  (t_3-2) (z_1^{-2} + z_2^{-2} + z_3^{-2}+ z_4^{-2})   - t_5 \Big) 
\eeq
and so on ...

\subsection{Specialisation to the Weil-Petersson measure}

Now, we specialize to the Weil-Petersson spectral curve of \cite{EOWP}:
\beq
y(z) = {1\over 2\pi}\sin{(2\pi z)}
\quad \to \quad
t_{2d+3} = {(2i\pi)^{2d}\over 2d+1!} + 2\delta_{d,0}
\quad \to \quad
f(z) = 1- \ee{-4\pi^2 z}
\eeq
so that:
\beq
\td{f}(z) = 4\pi^2 z 
\quad \to \quad
\td{t}_b = 4\pi^2 \,\delta_{b,1} + \delta_{b,0}\,\ln{(-2)}
\eeq
therefore each $b_i$ must be $1$, and we must have $k=d_0$, and we get:
\beq
\encadremath{
 V_{g,n}(P_1,\dots,P_n) 
 = 2^{-d_{g,n}} \!\!\!\!  \sum_{d_0+d_1+\dots+d_F=d_{g,n}}\, {2^{d_0}\over d_0!}\, \prod_F {P_F^{2 d_F}\over d_F !} \,\,< (2\pi^2\kappa_1)^{d_0}\, \prod_F  \psi_F^{d_F}>   
} \eeq
which is, after Wolpert's relation \cite{Wolpert}, the Weil-Petersson volume since $2\pi^2 \kappa_1$ is the Weil-Petersson K\" ahler form,
and thus, we have rederived Mirzakhani's recursion relation.

\subsection{Specialisation to the $\kappa_2$ measure}

Just to illustrate our method, we consider the integrals with only $\kappa_2$:
\beq
 V_{g,n}(P_1,\dots,P_n) 
 = 2^{-d_{g,n}} \!\!\!\!  \sum_{2d_0+d_1+\dots+d_F=d_{g,n}}\, {1\over d_0!}\, \prod_F {P_F^{2 d_F}\over d_F !} \,\,< (\td{t}_2\kappa_2)^{d_0}\, \prod_F  \psi_F^{d_F}>   
\eeq
which correspond to the conjugated times
\beq
\td{f}(z) =  \td{t}_2 z^2
\quad \to \quad
f(z) = \sum_{k=1}^\infty {(-1)^k\,\td{t}_2^k\over k!}\, z^{2k}
\eeq
i.e. $t_3=3$, and
\beq
t_{4a+3} = 4 (-1)^a\, \td{t}_2^{a}\,\,{2a!\over a! (4a+1)!}\,\, - \delta_{a,0}
\eeq
The corresponding spectral curve is:
\beq
y(z) = -{z\over 2} + 2 \sum_{k=1}^\infty (-\td{t}_2)^k\,{2k!\over k! (4k+1)!}\, z^{4k+1}
\eeq
with that spectral curve, the volumes $V_{g,n}$ satisfy the recursion of theorem \ref{mainth}.

\subsection{Specialisation to discrete measure}

Let us consider the example where $\L=\l \,{\rm Id}$, which is particularly important because
\beq
Z =\int dM \,\, \ee{-N\Tr({M^3\over 3} - M(\l^2+{1\over \l}))} 
\propto
\int dM \,\, \ee{-{N\over T}\Tr({1\over 2}M^2 - {M^3\over 3})} 
\eeq
where
\beq
T=-{1\over 8}\, (\l^2+{1\over \l})^{-3/2}
\eeq
i.e. Kontsevich integral reduces to the usual cubic one-matrix model, which is known to count triangulated maps \cite{BIPZ}.

In that case we have:
\beq
t_k = \l^{-k}
\eeq
thus for $b\geq 1$:
\beq
\td{t}_b = 2^b\l^{-2b}\,\sum_{l=1}^b {1\over l} (1-2\l^3)^{-l}\,\sum_{a_1+\dots+a_l=b,a_i>0} \prod_i (2a_i+1)!!
\eeq
For instance we have:
\beq
V_{0,3}(L) = {1\over \l^{-3}-2}
\eeq
\beq
V_{1,1}(L) = {1\over 2}\,{1\over \l^{-3}-2}\, (L<\psi_1>_1 +\td{t}_1<\kappa_1>_1)
= {-1\over 8(2-\l^{-3})}\,\left( {L\over 6} + {\l^{-5}\over 2-\l^{-3}} \right)
\eeq
where $\td{t}_1 = 6\l^{-2}\,(1-2\l^3)^{-1}$.

It would be interesting to understand how this relates to the discrete Regge measure on the set triangulated maps.
In the case of triangulated maps, loop equations, i.e. the recursion equation eq.\ref{recWgninmainth} are known as Tutte's equations \cite{Tutte} which give a recursive manner to enumerate maps.
This shows how general  the recursion equation eq.\ref{recWgninmainth} is.

\section{Other properties}

From the general properties of the invariants of \cite{EOFg}, we immediately have the following properties:
\begin{itemize}

\item Integrability. The $F_g$'s satisfy Hirota equations for KdV hierarchy.
That property is well known and it motivated the first works on Witten-Kontsevich conjecture \cite{Konts}.

\item Virasoro. The invariants of \cite{EOFg} were initialy obtained in \cite{eynloop1mat,EOFg} from the loop equations, i.e. Virasoro constraints satisfied by $Z(\L)$.

\item Dilaton equation, we have:
\beq
W_{g,n}(z_1,\dots,z_n) = {1\over 2g+n-2}\,\Res_{z\to 0} \Phi(z) W_{g,n+1}(z_1,\dots,z_n,z)
\eeq
where $d\Phi=y dx$.

For the Weil-Petersson case, after Laplace transform this translates into \cite{EOWP}:
\beq
V_{g,n}(L_1,\dots,L_n)_{\rm WP} = {1\over 2g+n-2}\,\,\, {\partial\over \partial L_{n+1}}\, V_{g,n+1}(L_1,\dots,L_n,2i\pi)_{\rm WP}
\eeq

\item It was also found in \cite{EOFg} how all those quantities behave at singular points of the spectral curve, and thus obtain the so-called double scaling limit.

\item The invariants constructed in \cite{EOFg} have many other nice properties, and it would be interesting to explore their applications to algebraic geometry...

\end{itemize}

\section{Conclusion}

In this paper we have shown how powerful the loop equation method is, and that the structure of the recursion equation eq.\ref{recWgninmainth} (i.e. Virasoro or W-algebra constraints) is very universal. 

We have thus provided a new proof of Mirzakhani's relations, exploiting the numerous properties of the invariants introduced in \cite{EOFg}.
However, the construction of \cite{EOFg} is much more general than that of Mirzakhani, since it can be applied to any spectral curve and not only the Weil-Petersson curve $y={1\over 2\pi}\,\sin{(2\pi \sqrt{x})}$.
In other words, we have Mirzakhani-like recursions for other measures, and theorem.\ref{mainth} gives the relationship between a choice of $t_k$'s (i.e. a spectral curve) and a measure on moduli spaces.
Moreover, the recursion relations always imply integrability and Virasoro.

It would be interesting to understand what the algebraic invariants $W_{g,n}$ defined by the recursion relation of \cite{EOFg} compute for an arbitrary spectral curve, not necessarily hyperelliptical neither rational...

\section*{Acknowledgments}
We would like to thank M. Bertola, J. Hurtubise, M. Kontsevich, D. Korotkin and N. Orantin for useful and fruitful discussions on this subject.
This work is partly supported by the Enigma European network MRT-CT-2004-5652, by the ANR project G\'eom\'etrie et int\'egrabilit\'e en physique math\'ematique ANR-05-BLAN-0029-01, by the Enrage European network MRTN-CT-2004-005616, 
by the European Science Foundation through the Misgam program,
by the French and Japaneese governments through PAI Sakurav, by the Quebec government with the FQRNT.

\end{document}